\documentclass[reqno]{amsart}
\usepackage{a4wide}
\usepackage[all]{xy}
\newtheorem{theorem}{Theorem}
\newtheorem{lemma}[theorem]{Lemma}

\theoremstyle{definition}
\newtheorem{definition}{Definition}
\newtheorem{example}{Example}
\theoremstyle{remark}
\newtheorem*{remark}{Remark}
\newcommand{\cc}{\mathbf{C}}
\newcommand{\hl}{\mathcal{H}^L}
\newcommand{\wl}{W^L}
\newcommand{\wlzero}{\wl_0}
\newcommand{\mlhalf}{m_{L/2}}
\newcommand{\mzerolhalf}{m_{0,L/2}}
\newcommand{\h}{\mathcal{H}}
\newcommand{\kk}{\mathcal{K}}
\newcommand{\z}{\mathbf{Z}}
\newcommand{\q}{\mathbf{Q}}
\newcommand{\R}{\mathbf{R}}
\newcommand{\rd}{\mathbf{R}^d}
\newcommand{\qp}{\mathbf{Q}_p}

\newcommand{\ztwo}{\mathbf{Z}_2}
\newcommand{\ftwo}{\mathbf{F}_2}
\newcommand{\qpd}{\mathbf{Q}_p^d}
\newcommand{\zpd}{\mathbf{Z}_p^d}
\newcommand{\ztwod}{\mathbf{Z}_2^d}
\newcommand{\ftwod}{\ftwo^d}
\newcommand{\rdxrd}{\mathbf{R}^d\times\mathbf{R}^d}

\newcommand{\qpdxqpd}{\mathbf{Q}^d_p\times\mathbf{Q}^d_p}
\newcommand{\zpdxzpd}{\mathbf{Z}^d_p\times\mathbf{Z}^d_p}
\newcommand{\qtwodxqtwod}{\mathbf{Q}^d_2\times\mathbf{Q}^d_2}
\newcommand{\ztwodxztwod}{\mathbf{Z}^d_2\times\mathbf{Z}^d_2}

\begin{document}
\title[Irreducible Models]{Models for the Irreducible Representation of a Heisenberg Group}
\author{Trond Digernes}
\address{Department of Mathematical Sciences\\The Norwegian
University of Science and Technology\\7491 Trondheim\\Norway}
\email{digernes@math.ntnu.no}
\urladdr{\href{http://www.math.ntnu.no/~digernes}{http://www.math.ntnu.no/~digernes}}
\thanks{This research was supported by the
Norwegian Research Council} \author{V. S.
Varadarajan}\address{Department of Mathematics\\University of
California\\Los Angeles, CA 90095\\USA}\email{vsv@math.ucla.edu}
\urladdr{\href{http://www.math.ucla.edu/~vsv}{http://www.math.ucla.edu/~vsv}}
\begin{abstract}
In its most general formulation a quantum kinematical system is
described by a Heisenberg group; the "configuration space" in this
case corresponds to a maximal isotropic subgroup. We study
irreducible models for Heisenberg groups based on \emph{compact}
maximal isotropic subgroups. It is shown that if the Heisenberg
group is 2-regular, but the subgroup is not, the "vacuum sector"
of the irreducible representation exhibits a fermionic structure.
This will be the case, for instance, in a quantum mechanical model
based on the 2-adic numbers with a suitably chosen isotropic
subgroup.

The formulation in terms of Heisenberg groups allows a uniform
treatment of $p$-adic quantum systems for all primes $p$, and
includes the possibility of treating adelic systems.
\end{abstract}
\maketitle
\section{Introduction}\label{intro}
In Hermann Weyl's formulation of quantum kinematics\footnote{All
the concepts discussed in this section will be explained in the
subsequent sections.}, a system with $d$ degrees of freedom is
described by operators $U(a)$ and $V(b)$ on $L_2(\rd)$, defined by
\[
(U(a)f)(x)=e^{ia{\cdot}x}f(x),\quad (V(b)f)(x)=f(x+b),\quad f\in
L_2(\rd)
\]
where
\[
\quad a=(a_1, \dots ,a_d),\,  b=(b_1,\dots ,b_d)\in\rd
\]
and
\[
x{\cdot}y=x_1y_1+\dots +x_dy_d \,.\] If we set
\[ W(a,b)=e^{(i/2)a{\cdot}b}U(a)V(b)\qquad (a, b)\in
{\R}^d\oplus{\R}^d
\]
then $(a,b)\mapsto W(a,b)$ is an irreducible \emph{projective}
unitary representation of $\rdxrd$ with multiplier $m$ given by
\[
m\left((a,b),(a',b')\right)=e^{(i/2)(a'{\cdot}b-a{\cdot}b')}
\]
i.e.,
\[
W(a,b)W(a',b')=e^{(i/2)(a'{\cdot}b-a{\cdot}b')}W(a+a',b+b')
\]
The original operators $U(a)$ and $V(b)$ are recovered from
$W(a,b)$ through $U(a)=W(a,0)$ and $V(b)=W(0,b)$. Focussing on
(the exponential of) the position operator, $U(a)$, we notice that
it is obtained by restricting $W$ to the \emph{configuration
space} $\rd\times(0)$, which is a \emph{maximal isotropic
subgroup} of $\rdxrd$ with respect to the \emph{symplectic}
multiplier $m$.

All of the above makes sense in the context of a locally compact
abelian group $G$ with a \emph{symplectic} multiplier $m$, i.e., a
\emph{Heisenberg} group $(G,m)$. So, following Weyl, one can take
the position that quantum kinematics, in its most general form, is
described by an irreducible projective representation $W$ of a
Heisenberg group $(G,m)$.

In this general setting there is no obvious candidate for the
"configuration space", but one can adopt the view that it
corresponds to a choice of a maximal isotropic subgroup $L$ for
$m$. One can then study the "position operator" by restricting the
projective representation $W$ to $L$, thereby obtaining an
ordinary representation $U$ of $L$.

The objective of this article is to analyze the representation $U$
in the case where $m$ has a \emph{compact} maximal isotropic
subgroup. Of course, this situation never occurs in the
conventional case $G=\rdxrd$, but it does occur, for instance, if
we choose our phase space to be $G=\qpdxqpd$. In the generic case
(i.e., both $G$ and $L$ are 2-regular), $U$ decomposes into
one-dimensional sub-representations. However, if $G$ is 2-regular,
but $L$ is not, an interesting phenomenon occurs: When considered
on the so-called vacuum space of $W$, the lift of $U$ to $L/2$
gives rise to a set of operators exhibiting a \emph{fermionic
structure} (Theorem~\ref{fermionic-theorem}).

The results of this article extend results obtained by Vladimirov,
Volovich and Zelenov for the case $G=\qpdxqpd$ (see their book
\emph{"$p$-adic Analysis and Mathematical Physics"}
\cite[244--247]{VVZ94}). In particular, we obtain a version of
their Theorem~3 on page~247 which is valid also for $p=2$; in
fact, this is a particularly interesting case, since this is where
the phenomenon of a "fermionic structure" occurs\footnote{Although
not stated explicitly in \cite{VVZ94}, Theorem~3 on page~247
requires $p\neq2$.}. Our approach also brings out more clearly the
mechanisms which are at work here.

The paper is organized as follows: In Section~2 we remind the
reader of the basic facts concerning conventional quantum
mechanics and Weyl systems. Section~3 contains a discussion on
multipliers, bicharacters, isotropic subgroups and models for the
unique irreducible representation of a Heisenberg group. Finally,
in Section~4, we present and prove our main results.

\section{Preliminaries}\label{prelim}
Quantum kinematics is based on the well known Heisenberg
commutation rules (with $\hbar =1$)
\[
[p_j,q_k]=-i\delta _{jk}I\qquad 1\le j,k\le d\tag{$H_d$}
\]
where $q_j,p_j$ $(1\leq j\leq d)$ are the position and momentum
coordinates of a quantum system with $k$ degrees of freedom.
Almost at the same time as these were discovered, Weyl noticed
that they  are the infinitesimal version of commutation rules
between the
   \emph{unitary groups} generated by the $q_j,p_k$; these are known
as the \emph{Weyl commutation rules.} Let
\[
U(a)=e^{ia{\cdot }x},\qquad  V(b)=e^{ib{\cdot }p}
\]
where
\[
\quad a=(a_1, \dots ,a_d),\qquad  b=(b_1,\dots ,b_d)
\]
and
\[
x{\cdot}y=x_1y_1+\dots +x_dy_d
\]
The Weyl commutation rules are then given by
\[\label{weylrels}
 U(a)V(b)=e^{-ia{\cdot}b}V(b)U(a)\qquad a,b\in\R \tag{$W_d$}
\]
A pair of unitary representations $U,V$ of ${\R}^d$ in a Hilbert
space is called a \emph{Weyl system} if they satisfy
\eqref{weylrels}. These two types of commutation rules are
formally equivalent and so Weyl took the point of view that
\eqref{weylrels} describes the kinematics of quantum systems whose
configuration space is a \emph{real} affine space. However the
concept of a Weyl system is much deeper than being just an
equivalent way to formulate quantum kinematics  as envisioned by
Heisenberg.  The point is, as Weyl himself discovered, that
although it was originally defined for ${\R}^d$, the notion of a
Weyl system can be formulated in much greater generality. In the
first place, the Weyl commutation rules involve just ${\R}^d$ and
the duality of ${\R}^d$ with itself given by the pairing
\[
a,b\longmapsto e^{ia{\cdot}b}
\]
and so one can speak of a Weyl system whenever one has a pair of
unitary representations $U,V$ of  abelian groups $A,B$ with a
nondegenerate pairing
\[
\langle,\rangle: A\times B\longrightarrow T\qquad (T=\{z\in
{\mathbf{C}}\ |\ |z|=1\})
\]
satisfying
\[
U(a)V(b)=\langle a,b\rangle^{-1}V(b)U(a)\qquad a\in A, b\in B
\]
One can then identify $B$ with a subgroup of the dual group
$\widehat A$ of $A$ by means of this pairing and then interpret
the pair $(U,V)$ as describing the quantum kinematics of a system
whose configuration space is the Abelian group $A$, whose momenta
lie in $B$, and which are covariant under the group of
translations of $A$. Weyl treated only the cases when $A$ and $B$
are either finite are finite dimensional Lie groups; both of these
are included in the above scheme if we stipulate that $A$ and $B$
are to be locally compact, and the resulting theory is completely
adequate for the treatment of quantum systems which are finite
dimensional; for treating quantum systems with \emph{infinitely}
many degrees of freedom, such as a quantum field, one has to give
up the local compactness of $A$ and $B$, and this introduces many
technical complications. We will not discuss this case here.

Weyl's generalization  already allows us  to describe quantum
systems whose configuration space can be very arbitrary, in
particular if it is an affine space over a \emph{nonarchimedean
local field $K$.}  We can take $A$ to be a finite dimensional
vector space $V$ over such a field, $B$ to be its (linear) dual
$V^\ast $, and $(\ ,\ )$ to be the pairing
\[
(a,b)= \chi (\langle a,b\rangle)
\]
where $\langle\ ,\ \rangle$ is the natural $K$--valued pairing of
$V$ and $V^\ast $, and $\chi $ is a basic character of $K$;  and
then work with Weyl systems for the pair $V,V^\ast $. This is  our
basic point of view in this article.

Weyl carried out a second generalization of the concept of a Weyl
system. To motivate this we observe that the definition
\eqref{weylrels} emphasizes the configuration space and its dual,
and so is really associated to a particular splitting of the phase
space. To get a more invariant description one should try to
formulate the concept of a Weyl system directly on the phase
space. Weyl noticed that the map\footnote{See Section~\ref{bichar}
for an explanation of the factor $e^{(i/2)a{\cdot}b}$.}
\[ (a,b)\longmapsto W(a,
b)=e^{(i/2)a{\cdot}b}U(a)V(b)\qquad (a, b)\in
{\mathbf{P}}={\R}^d\oplus{\R}^d
\]
is a \emph{projective unitary representation} of the phase space
$\mathbf{P}$, with the property that the phase factors which
determine the departure of $W$ from being an ordinary
representation, encode the \emph{symplectic structure} of the
phase space in the sense that they are of the form
\[
e^{i\beta (x,y)}
\]
where $\beta$ is the natural symplectic form on the phase space.
He was able to show that an affine space admits a faithful
irreducible projective representation only when it is symplectic,
and further that in this case the phase factors of the projective
representation are of the above form and hence determined by the
symplectic structure. He was thus led to his final and decisive
formulation of quantum kinematics as described by a
\emph{projective unitary representation of an abelian group with a
symplectic structure.}

Abelian groups with a symplectic structure are called
\emph{Heisenberg groups.} Our point of view can thus be described
as follows. \emph{Quantum kinematics for systems with finitely
many degrees of freedom is described by a projective unitary
representation that is canonically associated to a locally compact
Heisenberg group.} By specializing the Heisenberg group to one
defined over a local field or an adele ring we can then retrieve
the cases that are of greater interest, namely quantum systems
over local fields and adele rings. For Weyl's work, see
\cite{wey31}.

\subsection{The Schr\"odinger model for \eqref{weylrels} and its
uniqueness}\label{schrodinger} The  \emph{Schr\"odinger model} for
\eqref{weylrels} consists of the Hilbert space ${\h}=L^2({\R}^d)$
with
\[\label{sch}
(W(a,b)f)(t)=e^{(i/2)a{\cdot}b}e^{ia{\cdot}t}f(t+b)\tag{Sch}
\]
which is the same as the requirement
\[
(U(a)f)(t)=e^{ia{\cdot}t}f(t),\qquad (V(b)f)(t)=f(t+b) \] The
system of operators $\{ W(a,b)\}$ is irreducible. The central
question of quantum kinematics is  whether the commutation rules
\eqref{weylrels} can be satisfied only by the Schr\"odinger model,
i.e., whether \emph{up to irreducibility and unitary equivalence
the Schr\"odinger model is the unique solution of
\eqref{weylrels}.} The Heisenberg--Weyl commutation rules clearly
form the content of matrix mechanics while the Schr\" odinger
model formulates quantum kinematics in the form that is known as
wave mechanics. Weyl, who was the first to formulate the
uniqueness  question as the uniqueness of the pair $(U,V)$ or $W$,
clearly understood that its affirmative solution is what is needed
to show \emph{the equivalence of wave and matrix mechanics.} He
could not prove the uniqueness except in the special case when $a$
and $b$ vary in a finite cyclic group; it was proved a bit later
by Stone and Von Neumann for ${\R}^d$. \begin{theorem}[Stone--Von
Neumann] Any Weyl system for ${\R}^d\oplus {\R}^d$ is a direct sum
of copies of the Schr\"odinger model.
\end{theorem}
For Von Neumann's original proof see \cite{vn31,sto30}. For a more
recent treatment see \cite{var96}.

\section{Projective
representations}\label{projective} We give here a brief discussion
of projective representations, multipliers, and Heisenberg groups.

Let $G$ be a separable locally compact group and $T$ the
multiplicative group of complex numbers of absolute value $1$. An
\emph{$m$-representation} of $G$ is a (Borel) map of $G$ into the
unitary group of a Hilbert space such that for some (necessarily
Borel) map $m$, called the \emph{multiplier} of $U$, of $G\times
G$ with values in $T$,
\[
U(x)U(y)=m(x,y)U(xy),\qquad U(1)=1
\]
If we change $U$ to $U'=aU$, where $a$ is a Borel map of $G$ into
$T$, then $U'$ is an $m'$-representation where
\[\label{mequiv}
m'(x,y)=m(x,y)\frac{a(x)a(y)}{a(xy)}\tag{$\ast$}
\]
The multiplier itself satisfies
\begin{align}
m(x,1)&=m(1,x)=1\\
m(xy,z)m(x,y)&=m(x,yz)m(y,z)
\end{align}
Any Borel function $m$ on $G\times G$ with values in $T$
satisfying the above relations is called a \emph{multiplier} for
$G$; two multipliers $m,m'$ related as in \eqref{mequiv} are
called \emph{equivalent.} If
\[
m(x,y)=\frac{a(x)a(y)}{a(xy)}
\]
$m$ is called \emph{trivial}; in this case, if $U$ is an
$m$--representation, $aU$ becomes an ordinary representation.
Under pointwise multiplication, the multipliers form a group, the
trivial ones form a subgroup, and the quotient group of
equivalence classes of multipliers is called the \emph{multiplier
group} (of $G$). It is denoted by $M(G)$.

\subsection{Projective representations viewed as
ordinary representations of central extensions}\label{central} In
the above definitions we have considered no continuity assumptions
on multipliers. In most situations the multipliers are at least
locally continuous near $(1,1)$, and for the case we are
interested in, namely for Weyl systems and Heisenberg groups, even
continuous. For simplicity we assume here that $m$ is a
\emph{continuous}\footnote{This construction can be carried out
also when $m$ is only assumed to be Borel. It leads to the Weil
topology, the unique locally compact group topology on $G_m$ which
generates the product Borel structure on $G\times T$ (see
\cite{var85} for details of this theory which goes back to
\cite{mac58}).} multiplier. Then
\[
G_m=G\times T
\]
becomes a separable locally compact group (in the product
topology) with the multiplication
\[
(x,s)(y,t)=(xy,st\cdot m(x,y))
\]
For an $m$--representation $U$, the map
\[
U_m : (x,t)\longmapsto tU(x)
\]
is an ordinary representation of $G_m$ with $U_m(1,t)=t1$. We call
any representation of $G_m$ that restricts to $t\longmapsto t1$ on
the group $1\times T$, a \emph{basic representation}. So $U_m$ is
a basic representation of $G_m$, and
\[
U\longleftrightarrow U_m
\]
is a bijective correspondence between $m$--representations of $G$
and basic ordinary representations of $G_m$. The subgroup $1\times
T$ (which we identify with $T$) lies in the center of $G_m$ and
$G_m$ is a \emph{central extension of} $G$ by $T$ described by the
exact sequence
\[
1\longrightarrow T\longrightarrow G_m\longrightarrow
G\longrightarrow 1
\]

\subsection{Alternating bicharacters as multipliers}\label{bichar}
\emph{From now on we shall work with a separable locally compact
Abelian group $G$ written additively.}  We shall take a closer
look at the structure of the group of multipliers of $G$. First we
consider bicharacters. A \emph{bicharacter} of $G$ is a continuous
map $b$ of $G\times G$ into $T$ which is a character in each
argument; it induces a morphism $\beta =\beta _b$ of $G$ into
$\widehat G$ in the usual manner:
\[
(x,\beta (y))=b(x,y)
\]
The functional equations satisfied by a bicharacter are
sufficiently strong that one can show that \emph{any}    map $b$
of $G\times G$ into $T$ which is a  Borel homomorphism in each
argument, is a bicharacter (see \cite[p. 281]{mac58}). A
bicharacter $b$ is called \emph{alternating} or
\emph{antisymmetric} if it satisfies
\begin{gather}\label{bich}
b(x,y)b(y,x)=1,\qquad  b(x,x)=1
\end{gather}
and  \emph{symplectic} if it is alternating and $\beta _b$ is an
isomorphism.

Any bicharacter  is a continuous multiplier. Let $m$ be \emph{any
multiplier,} continuous or not, and let us set
\[
m^\sim (x,y)=\frac{m(x,y)}{m(y,x)}
\]
Then
\[
m^\sim (x,y)m^\sim (y,x)=1,\qquad  m^\sim (x,x)=1
\]
and it is easily checked that $m^\sim $ is a homomorphism in each
argument; in fact,
\begin{align*} m^\sim (x,y+z)&=\frac{m(x,y+z)}{m(y+z,x)}\\
&=\frac{m(x+y,z)m(x,y)}{m(y,z)m(y+z,x)}\\
 &=\frac{m(x+y,z)m(x,y)}{m(y,z+x)m(z,x)}\\
&=\frac{m(x,y)}{m(y,x)}\cdot \frac{m(y+x,z)m(y,x)}{m(y,z+x)m(z,x)}\\
&=\frac{m(x,y)}{m(y,x)}\cdot\frac{m(y,x+z)m(x,z)}{m(y,z+x)m(z,x)}\\
 &=m^\sim (x,y)m^\sim (x,z)
\end{align*}
while the homomorphism property in the first argument follows from
the antisymmetry. So, as $m^\sim $ is obviously Borel,   \emph{it
is an  alternating bicharacter.}  The bicharacter $m^\sim $ arises
naturally because of the fact that for any $m$-representation $U$
we have the commutation rule
\[
U(x)U(y)U(x)^{-1}U(y)^{-1}=m^\sim (x,y)1
\]
The following definition will be needed only for the case $p=2$,
but we state it for a general (prime number) $p:$
\begin{definition} Let $p$ be a prime number. An Abelian group $G$
is said to be \emph{$p$-divisible} (resp.\ \emph{$p$-injective})
if the map
\[x\in G\mapsto px\in G
\]
is surjective (resp.\ injective); if the map is bijective, we say
that $G$ is \emph{$p$-regular}.
\end{definition}
As an easy verification shows, a locally compact Abelian group $G$
is $p$-divisible (resp.\ $p$-injective) if and only if its dual
group $\widehat{G}$ is $p$-injective (resp.\ $p$-divisible). We
write $x\mapsto x/p$ for the inverse map of $x\mapsto px$ (when it
exists); more generally, if $A$ is a subset of $G$, we write $A/p$
or $(1/p)A$ for the inverse image of $A$ under the map $x\mapsto
px$ of $G$ into $G$.

\begin{lemma}\label{sqrt}If $G$ is $2$-regular and $\beta$ is a
bicharacter of $G$, there is a unique bicharacter $\beta^{1/2}$
such that $(\beta^{1/2})^2=\beta$.
\end{lemma}
\begin{proof}Define
\[
\beta^{1/2}(a,b)=\beta(a/2,b/2)^2\qquad a,b\in G
\]
Then $\beta^{1/2}$ is a bicharacter and $(\beta^{1/2})^2=\beta$.\\
\emph{Uniqueness.} Assume $\beta_1,\beta_2$ are two bicharacters
such that $\beta_1^2=\beta_2^2$, and set $\gamma=\beta_1/\beta_2$.
Then $\gamma^2=1$, so $\gamma(x,y)=\pm1$ for all $x,y\in G$. But
since $G$ is 2-regular, $\gamma(x,y)=\gamma(x/2,y)^2=1$ for all
$x,y$. Thus $\beta_1=\beta_2$, and uniqueness follows.
\end{proof}

Let $\Lambda ^2(G)$ be the multiplicative group of alternating
bicharacters. Then the map
\[
m\longmapsto m^\sim
\]
is a homomorphism of the group of multipliers of $G$ into the
group $\Lambda^2(G)$. We now have the following important lemma:
\begin{lemma}\label{lem:bichar}\mbox{}
\begin{enumerate}
\item\label{i} If $m_i\ (i=1,2)$ are multipliers, then
\[
m_1\simeq m_2\Longleftrightarrow m_1^\sim =m_2^\sim
\]
In particular, a multiplier $m$ is trivial if and only if $m^\sim
=1$, i.e., $m$ is symmetric. If $m$ is in addition continuous, we
can find a continuous map $a$ of $G$ into $T$ such that
\[
m(x,y)=\frac{a(x+y)}{a(x)a(y)}
\]
Indeed, any Borel map $a\colon G\longrightarrow T$ such that
$m(x,y)=a(x+y)/a(x)a(y)$ is automatically continuous.
\item\label{ii} The map
\[
m\longmapsto m^\sim
\]
is a homomorphism whose kernel is the group of trivial
multipliers, and so induces   an injection of $M(G)$ into $\Lambda
^2(G)$. \item\label{iii} If $G$ is $2$-regular, the map
\emph{(\ref{ii})} induces an isomorphism between $M(G)$ and
$\Lambda ^2(G)$. In this case, if $m\in M(G)$ with image $n\in
\Lambda ^2(G)$, then $n^{1/2}$ is the unique alternating
bicharacter in the multiplier class of $m$.
\end{enumerate}
\end{lemma}
\begin{proof}
Indeed, if $m$ is trivial, it is symmetric and  $m^\sim =1$.
Conversely, if $m^\sim =1$, then $m$ is symmetric and so $G_m$ is
Abelian. But then the dual of the exact sequence
\[
1\longrightarrow T\longrightarrow G_m\longrightarrow
G\longrightarrow 1
\]
is the exact sequence
\[
1\longrightarrow \widehat G\longrightarrow \widehat
{G_m}\longrightarrow {\mathbf{Z}}\longrightarrow 1
\]
As ${\mathbf{Z}}$ is free, this exact sequence splits, and so the
first exact sequence also splits, showing that $T$ is a direct
summand of $G_m$. Hence there is a continuous homomorphism
\[
x\longmapsto (x,a(x))
\]
of $G$ into $G_m$. Writing out the condition that this is a
homomorphism we get
\[
 m(x,y)=\frac{a(x+y)}{a(x)a(y)}
\]
In general $a$ is not continuous as $G_m$ need not have the
product topology. But if $m$ is continuous, $G_m$ \emph{has the
product topology,} and so $a$ is continuous. Since any two choices
of $a$ related as above to $m$ differ by a character, it follows
that any $a$ such that $m(x,y)=a(x+y)/a(x)a(y)$ is necessarily
continuous. If $m_1^\sim =m_2^\sim$, then $(m_1/m_2)^\sim =1$ so
that $m_1/m_2\simeq 1$, or $m_1\simeq m_2$. This proves (\ref{i})
and (\ref{ii}).\\
To prove  (\ref{iii}) note that $n^{1/2}$ is alternating and
$(n^{1/2})^\sim =(n^{1/2})^2=n$. The rest follows from
Lemma~\ref{sqrt}
\end{proof}
\begin{example}In Section~\ref{prelim} we defined the Weyl map
\[ (a,b)\longmapsto W(a,
b)=e^{(i/2)a{\cdot}b}U(a)V(b)\qquad (a, b)\in
{\mathbf{P}}={\R}^d\oplus{\R}^d
\]
which is a projective representation with multiplier
\[
m\left((a,b),(a',b')\right)=e^{(i/2)(a'{\cdot}b-a{\cdot}b')}
\]
A more natural definition would perhaps have been
\[ (a,b)\longmapsto W'(a,
b)=U(a)V(b)\qquad (a, b)\in {\mathbf{P}}={\R}^d\oplus{\R}^d
\]
which is a projective representation with multiplier
\[
m'\left((a,b),(a',b')\right)=e^{a'{\cdot}b}
\]
However, $m'$ is neither alternating nor nondegenerate, but it is
equivalent to $m$, which has both of those properties (i.e., $m$
is symplectic). The function $c(a,b)=e^{(i/2)a\cdot b}$ implements
the equivalence between $m$ and $m'$, and $m$ is the unique
alternating bicharacter in the multiplier class of $m'$, the
existence of which was guaranteed by part (\ref{iii}) of the lemma
above.
\end{example}
\begin{example} Similarly, if $A$ and $B$ are two Abelian groups with
a nondegenerate pairing
\begin{align*}
\langle\,\:,\:\rangle\colon A\times B\longrightarrow T&& (T=\{z\in
{\mathbf{C}}\ |\ |z|=1\})
\end{align*}
a Weyl system for $(A,B)$ was defined in Section~\ref{prelim} as a
pair of unitary representations $(U,V)$ satisfying
\begin{align*}
U(a)V(b)&=\langle a,b\rangle^{-1}V(b)U(a)&& a\in A,\: b\in B\\
\intertext{The corresponding Weyl operator}
W'(a,b)&=U(a)V(b)&&a\in A,\:b\in B\\
\intertext{is a projective representation of $G=A\times B$ with
multiplier}
m'\left((a,b),(a',b')\right)&=\langle a',b\rangle&&a,a'\in A,\:b,b'\in B\\
\intertext{Again, this multiplier is neither alternating nor
nondegenerate. But if $A$, and hence also $B$, is 2-regular, $m'$
is equivalent to the symplectic bicharacter}
m\left((a,b),(a',b')\right)&=\frac{\langle
a',b\rangle^{1/2}}{\langle
a,b'\rangle^{1/2}}&&a,a'\in A,\:b,b'\in B\\
\intertext{via the function $c\colon(a,b)\in G\mapsto\langle
a,b\rangle^{1/2}\in T$ (here the notation $\langle
a,b\rangle^{1/2}$ refers to the square root of the bicharacter
$\langle a,b\rangle$ as discussed in Lemma~\ref{lem:bichar}). The
representation} W(a,b)&=\langle a,b\rangle^{1/2}U(a)V(b)&&(a,b)\in
G
\end{align*}
is projective with multiplier $m$; and again, $m$ is the unique
alternating bicharacter in the multiplier class of $m'$.
\end{example}

\subsection{Heisenberg groups}\label{heisenberg}
Let $G$ be abelian and $m$ a multiplier for $G$. $G_m$ is called a
\emph{Heisenberg group} (with respect to $m$) and $m$ a
\emph{Heisenberg multiplier} (for $G$) if $m^\sim $ is a
symplectic bicharacter of $G$. This definition is prompted by the
observation that if $G={\R}^d\times {\R}^d$ and $m$ is the
multiplier of the projective representation $W$ in the
Schr\"odinger model, then $m^\sim $ is symplectic. If $G$ is
$2$--regular we shall usually replace $m$ by the unique
alternating bicharacter within the class of $m$ and so assume that
$m$ itself is a symplectic bicharacter. The most general form of
the uniqueness of projective unitary representations of abelian
groups with Heisenberg multipliers  is then the following theorem
proved by Mackey \cite{mac49,mac58}; see also
\cite{vn31,sto30,mum91,var96}:
\begin{theorem}[Mackey]If $m$ is a Heisenberg multiplier for $G$,
every $m$--representation of $G$ is a direct sum of copies of a
unique irreducible $m$--representation.
\end{theorem}
Mackey's theorem is thus the climax of the following succession of
ideas of Heisenberg, Weyl, and Mackey himself:
\begin{itemize}
\item\emph{Weyl system} \item\emph{projective unitary
representation of $A\oplus \widehat A$} \item\emph{projective
unitary representation of abelian group    with a Heisenberg
multiplier} \item\emph{basic representation of Heisenberg groups}
\end{itemize}
When $G=A\times \widehat A$ and the projective representation is
\[W(a,a^\ast )=U(a)V(a^\ast )\] where $U,V$ are unitary
representations of $A,\widehat A$ respectively, the unique
irreducible representation has  the Schr\"odinger model. In the
general case of an arbitrary $G$ with Heisenberg $m$ the model is
a little more subtle to describe
%(see Section~\ref{model}).

We shall not give the proof of this theorem here. The interested
reader may consult the above references; for a detailed
discussion, see \cite{var96} where a proof is given for the case
$G=A\times \widehat A$.

\subsection{Isotropic subgroups}\label{isotrop}
Let $m$ be a bicharacter on $G$ (we assume $m$ is Borel, and hence
continuous). For each subset $A$ of $G$ we define the \emph{polar}
$A'_m$ of $A$ with respect to $m$ by
\[A'_m=\{x\in G;\text{$m(x,a)=1$ for all $a\in A$}\}\]
We clearly have $ A\subset B\Longrightarrow B'_m\subset A'_m$ and
$A\subset A''_m$, and thus
\begin{align*}
A\subset &A''_m= A^{(4)}_m=A^{(6)}_m=\cdots,\\
%\text{and}\qquad\qquad
&A'_m=A'''_m=A^{(5)}_m=A^{(7)}_m=\cdots\,,
\end{align*} where $A^{(n)}_m$ denotes the $n$-th polar.
Polars are closed subgroups (since we are assuming $m$ is
continuous). A subgroup $A$ is said to be \emph{isotropic} for $m$
if $m|_{A\times A}\equiv1$, which is the same as saying that
$A\subset A'_m$. Maximal isotropic subgroups exist by Zorn's
lemma. Since the closure of an isotropic subgroup is again
isotropic, maximal isotropic subgroups are closed. $A$ is maximal
if and only if $A=A'_m$.

If $m$ is alternating, then
$m^{\sim}(x,y)=\frac{m(x,y)}{m(y,x)}=m^2(x,y)=m(2x,y)=m(x,2y)$,
and so we get the following relation between $A'_{m}$ and
$A'_{m^{\sim}}:$
\begin{align*}
A'_{m^{\sim}}=\{x\in G;m(x,2a)=1,\forall a\in
A\}=(2A)_m'=(1/2)A_m'\,,
\end{align*}
where, as usual, $(1/2)A_m'$ denotes the inverse image of $A_m'$
under the mapping $x\to 2x$. In particular, if $A$ is 2-regular,
we have
\[A'_{m^{\sim}}=A_m'\qquad\text{($A$ 2-regular).}\]
So in this case, maximal isotropy with respect to $m$ is the same
as maximal isotropy with respect to $m^\sim$.

\subsection{Model for the irreducible
$m$-representation}\label{model} We shall now discuss the
structure of the unique basic irreducible representation of a
locally compact Heisenberg group. The point is that although the
representation is unique, one can have many different models. Each
model brings to the foreground certain aspects of the
representation which are particularly transparent in that model.

We proceed by analogy with the Schr\"odinger model when
$G={\R}^d\oplus {\R}^d$ and try to split $G$ as a direct sum of a
``configuration space'' and its dual. This is not always possible
in the most general case, but we can approximate to this
situation.

Let $G$ be a separable  locally compact abelian group with a
Heisenberg multiplier $m$. Choose a subgroup $A$ which is maximal
isotropic with respect to $m^\sim$. Since $m^\sim $ is $1$ on
$A\times A$, it follows that $m$ is symmetric on $A\times A$ and
hence, by Lemma~\ref{lem:bichar},
\[ m(a,b)=\frac{c(a+b)}{c(a)c(b)}\qquad
a,b\in A
\]
for some Borel map $c$ of $A$ into $T$, which is continuous if $m$
is; the function $c$ is unique up to multiplication by a character
of $A$. The irreducible model is the Hilbert space
${\h}={\h}_{m,A,c}$ of functions $f$ on $G$ such that
\begin{align*}
f(x+a)=m(x,a)c(a)^{-1}f(x)\qquad x\in G, a\in A\tag{i}\\
\int _{G/A}|f([x])|^2d[x]<\infty \qquad ([x]=x+A)\tag{ii}
\end{align*}
The irreducible $m$--representation is defined by
\[
(W_m(y))f(x)=m(x,y)f(x+y)
\]
\begin{remark}By (i), $|f(x)|$ is constant on the $A$--cosets and so
defines a function $|f([x])|$ on $G/A$; it is this function that
appears under the integral sign in (ii).\end{remark} For a proof
that this is the irreducible version, see the discussion in
\cite{mum91,var96}.

\subsubsection{Examples}\label{examples}
Suppose that $G$ is $2$--regular, $m$  is a symplectic
bicharacter, and that $A$ and $B$ are maximal isotropic subgroups
such that $G\simeq A\times B$. Then $m$ induces an isomorphism of
$B$ with $\widehat A$ and hence we are in the setting of a Weyl
system for$(A, \widehat A)$. The model described is then the
Schr\"odinger model.  In  the general case the maximal isotropic
subgroup $A$ need not be a direct summand of $G$. This is the case
if
\[
G=V={\R}^d\times {\R}^d,\  m=e^{2\pi ib},\ A={\z}^d\times {\z}^d
\]
where $b$ is the skewsymmetric bilinear form with matrix
$\begin{pmatrix} 0&-I\\ I&0\end{pmatrix}$. In this case the model
leads to the classical theta functions \cite{mum91}. If
\[
G=V={\q}_p^d\times {\q}_p^d, \ m=\chi _p\circ b,\ A={\z}_p^d\times
{\z}_p^d
\]
where $b$ has the matrix above, then $A$ is maximal isotropic and
not a direct summand: it is a compact open subgroup, and $A=2A$
iff $p\not=2$. In the general adelic case  $G=V\times H$ where $V$
is a finite dimensional symplectic vector space over ${\R}$, and
$H$ is a totally disconnected group with a Heisenberg multiplier
$m$ admitting \emph{compact open subgroups which are maximal
isotropic for $m^\sim $.}

\section{Structure of vacuum sector}\label{vac}
In this section we perform a more detailed analysis of the unique
irreducible $m$-representation $W$. More specifically, we study
the case when the Heisenberg group $(G,m)$ has a maximal isotropic
subgroup $A$ which is \emph{compact}. The restriction of $W$ to
$A$ is an ordinary representation, and we analyze the so-called
\emph{vacuum sector} of this restriction. It turns out that the
results depend very much on the properties of the subgroup $A$: If
$A=2A$ (in particular, if $A$ is 2-regular), then there is a
unique (up to scalars) vacuum vector. On the other hand, if $A\neq
2A$, the vacuum sector has dimension larger than 1; in addition,
it has a canonical \emph{fermionic structure}.

Before discussing the two cases separately, we shall look into the
structure of $W$ a little more closely.

\subsection{Restriction to a compact maximal isotropic subgroup}
\label{restrict} Let $G$ be $2$--regular; $m$ a symplectic
bicharacter  for $G$; $L$ a maximal isotropic subgroup  for $m$
which is \emph{compact.} Then we have an isomorphism of $\widehat
L$ with $G/L$ such that $[y]:=y+L$ corresponds to the character
$a\longmapsto m(a,y)$ of $L$. Since $\widehat L$ is discrete, $L$
is \emph{open} in $G$. Let $W$ be an $m$--representation,
\emph{not necessarily irreducible.} Then the restriction of $W$ to
$L$ is an ordinary representation and so, writing
\[ {\h}_{[y]}=\{\psi :
W(a)\psi =m(a,y)\psi \ \forall a\in L\}
\]
we have
\begin{equation}\label{decomp}
{\h}=\oplus _{[y]\in G/L}{\h}_{[y]}
\end{equation}
Since
\[
W(a)W(x)\psi =m(a,2x)W(x)W(a)\psi=m(a,y+2x)W(x)\psi \qquad
\psi\in{\h}_{[y]},\: a\in L
\]
we see that
\begin{equation}\label{permute}
W(x) : {\h}_{[y]}\simeq {\h}_{[y+2x]}\qquad (x\in G)
\end{equation}
By \eqref{decomp} and \eqref{permute} ${\h}$ may be viewed as a
vector bundle on $G/L$ with a $G$-action   permuting the fibers
${\h}_{[y]}$ above the transitive action $x,[y]\longmapsto [y+2x]$
on $G/L$. In particular
\[
W\big|_L\simeq \dim ({\h}^L) \text{ times the regular
representation of } L
\]
where ${\h}^L={\h}_{[0]}$ is the subspace of vectors fixed by $L$.
\begin{lemma}\label{normalizer}
$L/2$ is the normalizer (in $G$) of $\hl$. Moreover, the action of
$W|_{L/2}$ on $\hl$ factors through $(L/2)/(2L)$.
\end{lemma}
\begin{proof}
Assume $x\in G$ normalizes $\hl$, i.e., $W(x)\hl\subset\hl$, and
let $0\neq\psi\in\hl$. Then, for $a\in L$,
\[
W(a)W(x)\psi=m(a,2x)W(x)W(a)\psi=m(a,2x)W(x)\psi
\]
Since $W(x)\psi\in\hl$, it follows that $m(a,2x)=1$. This holding
for all $a\in L$, and $L$ being maximal $m$-isotropic, we get
$2x\in L$, i.e., $x\in L/2$. On the other hand, if $x\in L/2$ and
$\psi\in\hl$, then
\[
W(a)W(x)\psi=m(a,2x)W(x)W(a)\psi=W(x)\psi
\]
for all $a\in L$; thus $W(x)\psi\in\hl$, i.e., $x$ normalizes
$\hl$.\\
As for the second statement, if $x\in L/2,a\in2L$, we have on
$\hl$:
\[W(x+a)=m(x,a)W(x)W(a)=m(x,a)W(x)=m(2x,a/2)W(x)=W(x)
\]
\end{proof}
\begin{remark}\label{remark_normalizer} It is because of the appearance of $L/2$ in this lemma
that our further analysis will depend very much on whether $L=2L$
or $L\neq2L$. Also, note that $(L/2)/(2L)$ is a \emph{finite}
group, since it is discrete ($2L$ being open in $L/2$) and compact
($L/2$ being compact). -In the sequel we will need to consider the
action of $W(L/2)$ on $\hl$ several times, and we introduce the
notation $\wl$ for this action; i.e., $\wl$ is shorthand for the
$m$-representation $W|_{L/2}$ acting on $\hl$.
\end{remark}

\begin{lemma}\label{decomp-invar}
For each closed $\wl$-invariant subspace $\kk\subset\hl$ let
$\h(\kk)$ be the smallest $W$-invariant, closed subspace in $\h$
containing $\kk$. Further, let $P\colon\h\to\hl$ be the orthogonal
projection. Then
\[
P\h(\kk)=\kk
\]
If $\hl$ is the orthogonal direct sum $\oplus_i\kk_i$ of
$\wl$-invariant closed subspaces $\kk_i$, then
\[
\h=\oplus_i\h(\kk_i)
\]
\end{lemma}
\begin{proof}
Let $\mathcal{L}$ be the linear span of the vectors
$\{W(x)\psi:x\in G,\psi\in\kk\}$. As
$W(y)W(x)\psi=m(y,x)W(x+y)\psi$, $\mathcal{L}$ is stable under
$W(G)$, and so $\overline{\mathcal{L}}=\h(\kk)$. Further, as
$W(x)\psi\in\h_{[2x]}$, we have
\[
PW(x)\psi=\begin{cases}0&\text{if } [2x]\neq0\\
W(x)\psi\in\kk&\text{if } [2x]=0
\end{cases}
\]
So $P\mathcal{L}\subset\kk$, and hence $P\h(\kk)\subset\kk$. As
the reverse inclusion is obvious,
the first statement follows.\\
Let $\kk_i$, $i=1,2$, be two orthogonal $W(L/2)$-invariant closed
subspaces of $\hl$. We claim that $\h(\kk_i)$ are orthogonal.  If
$y_i\in G$, $\psi _i\in {\kk}_i$ and $[2y_1]\neq[2y_2]$, then
$W(y_1)\psi _1$ and $W(y_2)\psi _2$ are orthogonal because they
lie respectively in ${\h}_{[2y_1]}$ and ${\h}_{[2y_2]}$, which are
orthogonal. If $2(y_1-y_2)\in L$, then  we have
\[
W(y_2)\psi _2=W(y_1)\psi _2' \qquad \psi
_2'=m(y_1,y_2-y_1)^{-1}W(y_2-y_1)\psi _2
\]
where $\psi _2'\in {\kk}_2$ (because $y_2-y_1\in L/2$), and hence
orthogonal to $\psi _1$. But then
\[
(W(y_2)\psi_2,W(y_1)\psi_1)=(W(y_1)\psi_2',W(y_1)\psi_1)=(\psi_2',\psi_1)=0
\]
Thus $W(y_1)\psi _1$ and $W(y_2)\psi _2$ are orthogonal for all
$y_i\in G$, $\psi _i\in {\kk}_i$, and so
${\h}({\kk}_1)\perp{\h}({\kk}_2)$. Since ${\h}$ is generated by
$\hl$ under $W$, the second statement is now clear.
\end{proof}
\begin{lemma}\label{decomp-irred}
$\hl$ is an orthogonal direct sum of subspaces $\kk_i$ which are
irreducible under $\wl$. $W$ is irreducible on $\h$ if and only if
$\hl$ is irreducible under the action of $\wl$. In particular, for
irreducible $W$,
\[
\dim(\hl)<\infty
\]
\end{lemma}
\begin{proof}
The first statement follows from the compactness of $L/2$ (or,
even more so, from the finiteness of $(L/2)/(2L)$), as does
the third (after proving the second statement).\\
As for the second statement: If $\hl$ is reducible under $\wl$,
then, by the previous lemma, $\h$ is reducible under $W$.
Conversely, if $Q$ is a nontrivial projection operator on $\h$
which commutes with $W(G)$, then $\hl=Q\hl\oplus(I-Q)\hl$ is a
nontrivial splitting of $\hl$ into $\wl$-invariant subspaces (if
the splitting were trivial, say $Q\hl=\{0\}$, we would have
$\{0\}=W(G)Q\hl=QW(G)\hl=Q\h$, a contradiction). This completes
the proof of the second statement and the lemma.
\end{proof}

\subsection{Vacuum vectors and coherent states when
{\mathversion{bold}$L=2L$}}\label{2reg} Since now $L/2=L$, $\wl$
acts trivially on $\hl$, and thus is irreducible if and only if
$\dim(\hl)=1$. Lemma~\ref{decomp-irred} gives
\begin{theorem}\label{thm:2reg} If $L=2L$, $W$ acts irreducibly on $\h$ if and
only if $\dim(\hl)=1$.
\end{theorem}
The state defined by the one--dimensional space ${\h}^L$ for
irreducible $W$ may be called the \emph{vacuum state}; the states
defined by the one--dimensional spaces ${\h}_{[y]}$ are called the
\emph{coherent states}; and the decomposition \eqref{decomp} gives
the \emph{coherent state structure} of $W$.

\subsection{Structure of the vacuum space when
{\mathversion{bold}$L\neq2L$}} \label{non2reg} We now assume that
$W$ is an \emph{irreducible} $m$-representation in the Hilbert
space $\h$; it follows that $\wl$ is irreducible
(Lemma~\ref{decomp-irred}). We write $\mlhalf$ for the restriction
of the multiplier $m$ of $W$ to $L/2\times L/2$.
\begin{lemma}\label{m-lift} The bicharacter $\mlhalf^\sim$ is the
lift of a symplectic bicharacter of $(L/2)/L$.
\end{lemma}
\begin{proof}
For $x,x'\in L/2$, $a,a'\in L$ we have
\begin{align*}\mlhalf^\sim(x+a,x'+a')&=\mlhalf^\sim(x,x')\mlhalf^\sim(x,a')\mlhalf^\sim(a,x')
\mlhalf^\sim(a,a')\\
&=\mlhalf^\sim(x,x')m(2x,a')m(a,2x')m(a,a')^2=\mlhalf^\sim(x,x')\,,
\end{align*}
showing the lifting property. Let $n$ denote the corresponding
bicharacter on $(L/2)/L$. If $x\in L/2$ and $n(x+L,y+L)=1$ for all
$y\in L/2$, then
%\begin{align*}
$m(x,2y)=n(x+L,y+L)=1$ for all $y\in L/2$, hence $x\in L$ by the
maximal $m$-isotropy of $L$. So $n$ is non-degenerate, hence
symplectic (since $(L/2)/L$ is finite).
%\end{align*}
\end{proof}
\begin{lemma}\label{lem:factoring}There is a projective representation
$\wl_0$ of $(L/2)/L$ on $\hl$ whose lift to $L/2$ is projectively
equivalent to $\wl$. The multiplier $\mzerolhalf$ of $\wl_0$ is
characterized up to equivalence  by the condition that the lift of
$\mzerolhalf^\sim$ to $L/2\times L/2$ is equal to $\mlhalf^\sim $.
In particular, $\mzerolhalf$ is Heisenberg for the group
$(L/2)/L$.
\end{lemma}
\begin{proof}
If $x\in L/2$ and $a\in L$, we have on $\hl$
\[ W(x+a)=m(a,x)W(x)
\]
So the morphism defined by $\wl$ of $L/2$ into the projective
group of $\hl$ factors through $L$. As $(L/2)/L$ is finite, it is
clear that we can choose a projective representation of $(L/2)/L$
on $\hl$ which corresponds to the above morphism of $(L/2)/L$.
Clearly the lift of $m_{0,L/2}^\sim$ must coincide with
$\mlhalf^\sim$. The preceding lemma shows that $\mzerolhalf$ is
Heisenberg.
\end{proof}
For an illustration of Lemma~\ref{lem:factoring}, and its proof,
see Figure~\ref{fig:diagr}.
\begin{figure}\centering
\parbox{4cm}{\xymatrix{L/2\ar[r]^\wl\ar[d]_\pi\ar[dr]^{\overline{\wl}}&\mathcal{U}(\hl)\ar[d]^\psi\\
(L/2)/L\ar[r]_{\overline{\wl_0}}&\overline{\mathcal{U}(\hl)}}}
\hspace{4cm}
\parbox{4cm}{\xymatrix{&\mathcal{U}(\hl)\ar[d]^\psi\\
(L/2)/L\ar[ur]^{\wl_0}\ar[r]_{\overline{\wl_0}}&\overline{\mathcal{U}(\hl)}}}
\caption[Illustration of
Lemma~\ref{lem:factoring}]{\emph{Illustration of
Lemma~\ref{lem:factoring}.} $\mathcal{U}(\hl)$ is the unitary
group of $\hl$,
$\overline{\mathcal{U}(\hl)}=\mathcal{U}(\hl)/\mathbf{T}$ the
corresponding projective group, $\pi$ and $\psi$ are the canonical
projection maps, and the remaining symbols are self-explanatory.
The relation $\psi\circ\wl_0\circ\pi=\psi\circ\wl$ forces the
existence of a function $a\colon L/2\to\mathbf{T}$ such that
$\wl_0\circ\pi(x)=a(x)\wl(x)$, $x\in L/2$, which in turn forces
the relation
$\mzerolhalf(\pi(x),\pi(y))=\frac{a(x)a(y)}{a(x+y)}\mlhalf(x,y)$,
giving $\mzerolhalf^\sim(\pi(x),\pi(y))=\mlhalf^\sim(x,y)$,
$x,y\in L/2$. }\label{fig:diagr}
\end{figure}
\begin{theorem}\label{thm-carrierspace}The vacuum space $\hl$ carries the unique
irreducible $\mzerolhalf$-representation of $(L/2)/L$. Moreover
\[
\dim(\hl)=2^d\qquad |(L/2)/L|=2^{2d}
\]
for some positive integer $d$.
\end{theorem}
\begin{proof}
The first statement follows at once from the lemmas above. Now
$V_2:=(L/2)/L$ is a finite dimensional vector space over the field
$\mathbf{F}_2:=\z/2\z$ (since all elements of $V_2$ have order 2),
and $\mzerolhalf^\sim$ has values $\pm1$; so we can view it as a
\emph{symplectic form} on $V_2$. Thus $V_2$ has even dimension
$2d$ over $\mathbf{F}_2$, and hence cardinality $2^{2d}$. Let $A$
be a maximal $\mzerolhalf^\sim$-isotropic subspace of $V_2$. Then
$A$ has dimension $d$ and cardinality $2^d$. The unique
irreducible $\mzerolhalf^\sim$-representation of $V_2$ can be
realized on the Hilbert space $\ell_2(A)$, which has complex
dimension $2^d$. Thus $\dim(\hl)=\dim(\ell_2(A))=2^d$.
\end{proof}

The following theorem shows that, in the case of $L\neq2L$, the
vacuum space comes equipped with a canonical fermionic structure.
\begin{theorem}[Fermionic structure]\label{fermionic-theorem}
For the projective representation $\wl_0$ of $(L/2)/L$ there is a
basis
 $(e_i)_{i=1}^{2d}$ for $(L/2)/L$ such that:
\begin{equation}\label{fermionic-relations}\wlzero(e_i)^2=1,\qquad \wlzero(e_i)\wlzero(e_j)
=-\wlzero(e_j)\wlzero(e_i),\:i\neq j
\end{equation}
\end{theorem}
\begin{remark}The relations (\ref{fermionic-relations}) define the basic
representation of the Clifford algebra in $2d$ dimensions in terms
of $2d$ \emph{anticommuting spin observables}. This gives the
\emph{fermionic structure} on $\hl$.
\end{remark}
To prove the theorem we need the following lemma:
\begin{lemma}Let $V$ be a finite dimensional vector space over the field
$\mathbf{F}_2$, and let $b$ be a symplectic form on $V$. Then
there is a basis $(e_i)$ for $V$ such that
\[b(e_i,e_j)=1-\delta_{ij}
\]
Moreover, $b=b_1-b_1^t$ where $b_1^t$ is the transpose of $b_1$
and $b_1$ is the bilinear form given by
\[b_1(e_i,e_j)=\begin{cases}0,&i\leq j\\1,&i>j\end{cases}
\]
\end{lemma}
\begin{proof}We first check that the bilinear form $b_0$ on
$V\times V$ given by
\[b_0(f_i,f_j)=1-\delta_{ij}\,,
\]
where $(f_i)$ is a basis, is symplectic. For this it suffices to
prove that $b_0$ is nondegenerate. Suppose $v=\sum_{1\leq
i\leq2d}a_if_i$ is orthogonal to all of $V$, and set
$s=\sum_{1\leq i\leq2d}a_i$. Then
\[b_0(v,f_j)=\sum_{i\neq j}a_i=s-a_j=0\quad\forall j\,,
\]
so $a_j=s=2da_1=0$, $\forall j$. Since any two symplectic forms
are isomorphic, the first statement is clear. The second is an
explicit calculation.
\end{proof}
\begin{proof}[Proof of Theorem~\ref{fermionic-theorem}]As noted above,
the bicharacter $\mzerolhalf^\sim$ gives rise to a unique
symplectic form $b$ on $V_2$, the relation between them being
$\mzerolhalf^\sim=\chi\circ b$, where $\chi$ is the non-trivial
character on $\ftwo$. With $b_1$ as in the previous lemma we have
$\mzerolhalf^\sim=(\chi\circ b_1)^\sim$. Hence $\mzerolhalf$ and
$\chi\circ b_1$ are equivalent so that we may assume
$\mzerolhalf=\chi\circ b_1$. If the basis $(e_i)$ is chosen as in
the previous lemma, we have
\begin{align*}\wlzero(e_i)^2&=\mzerolhalf(e_i,e_i)\wlzero(2e_i)=1\\
\wlzero(e_i)\wlzero(e_j)&=\mzerolhalf^\sim(e_i,e_j)\wlzero(e_j)\wlzero(e_i)=\chi(b(e_i,e_j))\wlzero(e_j)\wlzero(e_i)\\
&=\chi(1)\wlzero(e_j)\wlzero(e_i)=-\wlzero(e_j)\wlzero(e_i)\,,\quad
i\neq j
\end{align*}
\end{proof}

\section{Examples involving $\qp$}\label{examples2}
Set
\[G=\qpdxqpd,\quad
m(x,y)=\chi_p(x_1\cdot y_2-x_2\cdot y_1)
\]
where $x=(x_1,x_2)$ and $y=(y_1,y_2)$ are in $G$ and $\chi_p$ is a
basic character on $\qp$. The representation $W$ given by
\[(W(y)f)(s)=\chi_p(2s\cdot y_2+y_1\cdot y_2)f(s+y_1),\quad
y=(y_1,y_2)\in G,\,s\in\qpd,\,f\in L_2(\qpd)
\]
is an $m$-representation of $G$ on $\h=L_2(\qpd)$, and
\[L=\zpdxzpd\]
is compact subgroup of $G$ which is maximal isotropic for $m$. As
for the vacuum space $\hl$ we have $f\in \hl\iff W(a_1,a_2)f=f$
for all $(a_1,a_2)\in L$. Written out this becomes
\begin{align}
\chi_p(2s\cdot a_2+a_1\cdot a_2)f(s+a_1)&=f(s)\notag\\
f(s+a_1)&=\chi_p(2s\cdot a_2)^{-1}f(s),\quad
s\in\qpd,\,(a_1,a_2)\in L\label{cond:hl}
\end{align}
Setting $a_1=0$ and assuming $f(s)\neq0$, this gives
$\chi_p(2sa_2)=1$ for all $a_2\in\zpd$, and thus $s\in\zpd/2$;
i.e., all functions in $\hl$ have support in $\zpd/2$. Further,
setting $a_2=0$ in \eqref{cond:hl} and letting $a_1$ be arbitrary,
we get
\[f(s+a_1)=f(s)\quad s\in\zpd/2,\;a_1\in\zpd
\]
All in all we have:
\begin{align}
\hl&=\{f\in L_2(\zpd/2)\colon
f(s+a_1)=f(s),s\in\zpd/2,a_1\in\zpd\}\\
&\simeq l_2\left((\zpd/2)/\zpd\right)
\end{align}

\subsection{$G=\qpdxqpd$ with $p\neq2$}\label{pnot2}
In this case $\zpd/2=\zpd$ and so $(\zpd/2)/\zpd$ reduces to a
point, i.e., $\hl\simeq\cc$, confirming the result of
Theorem~\ref{thm:2reg}.

\subsection{$G=\qtwodxqtwod$}\label{p2}
Since $(\ztwo/2)/\ztwo\simeq\ftwo$ (the field of two elements),
$\hl$ is now isomorphic to the space of functions on $\ftwod$. To
identify the representation $\wlzero$ of Lemma~\ref{lem:factoring}
with the fermionic properties \eqref{fermionic-relations},
consider the action $W'$ of $(\ztwod/2)\times(\ztwod/2)$ on $\hl$
defined as $W'(x_1,x_2)=W(x_1,0)W(0,x_2)$. Written out this
becomes
\begin{gather}
(W'(x_1,x_2)f)(s)=\chi_2(2s\cdot x_2+2x_1\cdot x_2)f(s+x_1)\qquad
x_1,x_2,s\in\ztwod/2,\quad f\in \hl
\end{gather}
$W'$ is an irreducible projective representation with multiplier
$m'\big((x_1,x_2),(y_1,y_2)\big)=\chi_2(-2y_1\cdot x_2) $. Since
$m'^\sim\big((x_1,x_2),(y_1,y_2)\big)=\chi_2(2x_1y_2-2y_1\cdot
x_2)=m^\sim\big((x_1,x_2),(y_1,y_2)\big)$, $W'$ and $W$ are
projectively equivalent. $W'$ factors through $\ztwodxztwod$, and
the corresponding projective representation $W''$ of
$\ftwod\times\ftwod$ on $\hl\simeq l_2(\ftwod)$ is given by:
\begin{align*}
\big(W''(\pi(x_1),\pi(x_2))\tilde{f}\big)\big(\pi(s)\big)&=
(W'(x_1,x_2)f)(s)=\chi_2(2s\cdot x_2+2x_1\cdot x_2)f(s+x_1)\\
&=\chi(\pi(s)\cdot
\pi(x_2)+\pi(x_1)\cdot\pi(x_2))\tilde{f}(\pi(s)+\pi(x_1))
\end{align*}
where $\pi\colon\ztwo/2\to(\ztwo/2)/\ztwo=\ftwo$ is the canonical
map, $f=\tilde{f}\circ\pi$, and $\chi$ is the nontrivial character
on $\ftwo$. For the multiplier $m''$ of $W''$ we have:
\begin{align*}
m''\big((\pi(x_1),\pi(x_2)),(\pi(y_1),\pi(y_2))\big)
&=m'\big((x_1,x_2),(y_1,y_2)\big)=\chi_2(-2y_1\cdot x_2)
=\chi\big(\pi(y_1)\cdot\pi(x_2)\big)\\
m''^\sim\big((\pi(x_1),\pi(x_2)),(\pi(y_1),\pi(y_2))\big)
&=\chi\big(\pi(y_1)\cdot\pi(x_2)-\pi(y_2)\cdot\pi(x_1)\big)
=\chi_2(2y_2\cdot x_1-2y_1\cdot
x_2)\\
&=m'^\sim\big((x_1,x_2),(y_1,y_2)\big)=m^\sim\big((x_1,x_2),(y_1,y_2)\big)
\end{align*}
Working directly on $\ftwod$, the above formulae simplify to:
\begin{align*}
\left.
\begin{aligned}
(W''(a_1,a_2)f)(c)&=\chi(c\cdot a_2+a_1\cdot a_2)f(c+a_1)\\
m''((a_1,a_2),(b_1,b_2))&=\chi(b_1\cdot a_2)\\
m''^\sim((a_1,a_2),(b_1,b_2))&=\chi(b_1\cdot a_2-b_2\cdot a_1)
\end{aligned}\right\}\qquad\begin{aligned}
f&\in l_2(\ftwod)\\
a_i,b_i,c&\in\ftwod,\quad i=1,2
\end{aligned}
\end{align*}

\bibliographystyle{amsalpha}

\begin{thebibliography}{Mum91}

\bibitem[Mac49]{mac49}
George~W. Mackey, \emph{A {T}heorem of {S}tone and von {N}eumann},
Duke Math.
  J. \textbf{16} (1949), 313--326.

\bibitem[Mac58]{mac58}
\bysame, \emph{Unitary representations of group extensions, {I}},
Acta Math.
  \textbf{99} (1958), 265--311.

\bibitem[Mum91]{mum91}
David Mumford, \emph{Tata lectures on theta. {III}}, Progress in
Mathematics,
  vol.~97, Birkh\"auser Boston Inc., Boston, MA, 1991, With the collaboration
  of Madhav Nori and Peter Norman. \MR{93d:14065}

\bibitem[Sto30]{sto30}
Marshall~H. Stone, \emph{Linear {T}ransformations in {H}ilbert
{S}pace, iii},
  Proc. Nat. Acad. Sci. U.S.A. \textbf{16} (1930), 172--175.

\bibitem[Var85]{var85}
V.~S. Varadarajan, \emph{Geometry of quantum theory}, second ed.,
  Springer-Verlag, New York, 1985. \MR{87a:81009}

\bibitem[Var96]{var96}
\bysame, \emph{Quantum kinematics and projective unitary
representations of
  abelian groups}, Analysis, geometry and probability, Texts Read. Math.,
  vol.~10, Hindustan Book Agency, Delhi, 1996, pp.~362--396. \MR{98j:81114}

\bibitem[vN31]{vn31}
John von Neumann, \emph{Die {E}indeutigkeit der {S}chr\"o
dingerschen
  {O}peratoren}, Math. Ann. \textbf{104} (1931), 570--578.

\bibitem[VVZ94]{VVZ94}
V.~S. Vladimirov, I.~V. Volovich, and E.~I. Zelenov,
\emph{$p$-adic analysis
  and mathematical physics}, World Scientific Publishing Co. Inc., River Edge,
  NJ, 1994. \MR{95k:11155}

\bibitem[Wei64]{wei64}
Andr{\'e} Weil, \emph{Sur certains groupes d'op\'erateurs
unitaires}, Acta
  Math. \textbf{111} (1964), 143--211. \MR{29 \#2324}

\bibitem[Wei65]{wei65}
\bysame, \emph{Sur la formule de {S}iegel dans la th\'eorie des
groupes
  classiques}, Acta Math. \textbf{113} (1965), 1--87. \MR{36 \#6421}

\bibitem[Wey31]{wey31}
H.~Weyl, \emph{Theory of groups and quantum mechanics}, Dover, New
York, 1931,
  Ch. III, \S 16, Ch. IV, \S\S 14, 15.

\end{thebibliography}
\providecommand{\bysame}{\leavevmode\hbox
to3em{\hrulefill}\thinspace}
\providecommand{\MR}{\relax\ifhmode\unskip\space\fi MR }
\providecommand{\MRhref}[2]{%
  \href{http://www.ams.org/mathscinet-getitem?mr=#1}{#2}
} \providecommand{\href}[2]{#2}

\end{document}